\documentclass[a4paper,12pt]{article}
\begin{document}

\title{Quasi-permutable normal operators in octonion Hilbert spaces and spectra}
\author{Ludkovsky S.V.}
\date{15 March 2012}
\maketitle

\begin{abstract}
Families of quasi-permutable normal operators in octonion Hilbert
spaces are investigated. Their spectra are studied. Multiparameter
semigroups of such operators are considered. A non-associative
analog of Stone's theorem is proved. \footnote{key words and
phrases: non-commutative functional analysis, hypercomplex numbers,
quaternion skew field, octonion algebra, quasi-permutable operators,
spectra, spectral measure, non-commutative integration
 \\
}
\end{abstract}

\section{Introduction}
The theory of bounded and unbounded normal operators over the
complex field is classical and have found many-sided applications in
functional analysis, differential and partial differential equations
and their applications in the sciences
\cite{danschw,jungexu10,kadring,killipsimon09,zelditch09}.
Nevertheless, hypercomplex analysis is fast developing, because it
is closely related with problems of theoretical and mathematical
physics and of partial differential equations
\cite{brdeso,gilmurr,guesprqa}. On the other hand, the octonion
algebra is the largest division real algebra in which the complex
field has non-central embeddings \cite{dickson,baez,kansol}. The
octonion algebra also is intensively used in mathematics and various
applications \cite{emch,guetze,girard,krausryan,kravchot}.
\par Previously analysis over quaternion and octonions was
developed and spectral theory of bounded normal operators and
unbounded self-adjoint operators was described \cite{ludoyst,ludfov,
lujmsalop,lufjmsrf,ludancdnb}. Their applications in partial
differential equations were outlined
\cite{luspraaca,ludspr,ludspr2,ludvhpde,ludcmft12}. This paper is
devoted to families of quasi-permutable normal operators in octonion
Hilbert spaces. Their spectra are studied. Multiparameter semigroups
of such operators are considered. A non-associative analog of
Stone's theorem is proved.
\par Notations and definitions of papers
\cite{ludoyst,ludfov,lujmsalop,lufjmsrf,ludancdnb} are used below.
The main results of this article are obtained for the first time.

\section{Quasi-permutability of normal operators}
\par {\bf 1. Definitions.}
If $\mbox{}_jA$ is a set of $\bf R$ homogeneous ${\cal A}_v$
additive operators with ${\cal A}_v$ vector domains ${\cal
D}(\mbox{}_jA)$ dense in a Hilbert space $X$ over the Cayley-Dickson
algebra ${\cal A}_v$, $~2\le v$, $~j\in \Lambda $, $~\Lambda $ is a
set, then we denote by $alg_{{\cal A}_v} (\mbox{}_jA: ~ j\in \Lambda
)$ a family of all operators $B$ with ${\cal A}_v$ vector domains in
$X$ obtained from $(\mbox{}_jA: ~ j\in \Lambda )$ by a finite number
of operator addition, operator multiplication and left and right
multiplication of operators on Cayley-Dickson numbers $b\in {\cal
A}_v$ or on $bI$, where $I$ denotes the unit operator on $X$.

\par Let $\mbox{}_1A$ and $\mbox{}_2A$ be two normal operators in a
Hilbert space $X$ over the Cayley-Dickson algebra ${\cal A}_v$,
$~2\le v$. Suppose that $\mbox{}_1A$ and $\mbox{}_2A$ are affiliated
with a quasi-commutative von Neumann algebra $\sf A$ over ${\cal
A}_v$ with $~2\le v\le 3$. Let $\mbox{}_1E$ and $\mbox{}_2E$ be
their ${\cal A}_v$ graded projection valued measures defined on the
Borel $\sigma $-algebra of subsets in ${\cal A}_v$ (see also \S 2
and \S \S I.2.58 and I.2.73 in \cite{ludopalglamb}). In this section
the simplified notation $E$ instead of ${\hat {\bf E}}$ will be
used.
\par We shall say that two normal operators $\mbox{}_1A$ and $\mbox{}_2A$ quasi-permute if
$$(1)\quad \mbox{}_1E(\delta _1)\mbox{ }_2E(\delta _2) =\mbox{ }_2E(\delta _2)
\mbox{}_1E(\delta _1)$$ for each Borel subsets $\delta _1$ and
$\delta _2$ in ${\cal A}_v$.

\par Operators $A$, $\mbox{}_1A$ and $\mbox{}_2A$ are said to have
property $P$ if they satisfy the following four conditions
$(P1-P4)$:
\par $(P1)$ they are normal, \par $(P2)$ they are affiliated with a
von Neumann algebra $\sf A$ over either the quaternion skew field or
the octonion algebra ${\cal A}_v$ with $2\le v \le 3$ and
$$(P3)\quad A = \mbox{}_1A \mbox{ }_2A\mbox{  and}$$
$(P4)$ the family $alg_{{\cal A}_v} (I,A,A^*,\mbox{ }_1A,\mbox{
}_1A^*,\mbox{ }_2A,\mbox{ }_2A^*)=: {\sf Q}(A,\mbox{ }_1A,\mbox{
}_2A)=: \sf Q $ over ${\cal A}_v$ generated by these three operators
is quasi-commutative, that is a von Neumann algebra \par $cl
[alg_{{\cal A}_v} (I,AE(\delta ),A^*E(\delta ),\mbox{ }_1A\mbox{
}_1E(\delta _1),\mbox{ }_1A^*\mbox{ }_1E(\delta _1),\mbox{
}_2A\mbox{ }_2E(\delta _2),\mbox{ }_2A^*\mbox{ }_2E(\delta
_2))]\subset L_q(X)$
\\ contained in $L_q(X)$ is quasi-commutative for each bounded Borel
subsets $\delta $, $\delta _1, \delta _2\in {\cal B}({\cal A}_r)$,
where $2\le v\le 3$.
\par  It is possible to consider
a common domain ${\cal D}^{\infty }({\sf Q}) := \cap_{T\in {\sf Q}}
{\cal D}^{\infty }(T)$ for a family of operators $\sf Q$, where
${\cal D}^{\infty }(T) := \cap_{n=1}^{\infty } {\cal D}(T^n)$. Then
the family $\sf Q$ on ${\cal D}^{\infty }({\sf Q})$ can be
considered as an ${\cal A}_v$ vector space. Take the decomposition
${\sf Q}={\sf Q}_0i_0\oplus {\sf Q}_1i_1\oplus ... \oplus {\sf
Q}_{2^v-1}i_{2^v-1}$ of this ${\cal A}_v$ vector space with pairwise
isomorphic real vector spaces ${\sf Q}_0, {\sf Q}_1,...,{\sf
Q}_{2^v-1}$. Then as in \S 2.5 \cite{ludunbnormopla12} for each
operator $B\in \sf Q$ we put
$$(2)\quad B=\sum_j \mbox{}^jB\mbox{  with  }\mbox{}^jB={\hat \pi }^j(B)\in {\sf
Q}_ji_j$$ for each $j$, where ${\hat \pi }^j: {\sf Q}\to {\sf
Q}_ji_j$ is the natural $\bf R$ linear projection, real linear
spaces ${\sf Q}_ji_j$ and $i_j{\sf Q}_j$ are considered as
isomorphic, so that
$$(3)\quad \sum_{k=0}^{2^v-1} \mbox{ }^k{\hat T} = T.$$
\par If $E$ is an ${\cal A}_v$ graded projection valued measure on
the Borel $\sigma $-algebra ${\cal B}({\cal A}_v)$ for a normal
operator $T\in \sf Q$, for uniformity of this notation we put also
\par $(4)$ $\mbox{}^k{\hat E}(dz).ty={\hat \pi }^kE(dz).ty$ \\
for every vector $y\in X$ and $t=t_0i_0+...+t_{2^v-1}i_{2^v-1}\in
{\cal A}_v$, where $z\in {\cal A}_v$, $~t_0,...,t_{2^v-1}\in \bf R$,
$~E(dz).ty=E(dz).(ty)$.

\par {\bf 2. Lemma.} {\it Let operators $A$, $B$ and $D$ have property $P$
and let $F$ be an ${\cal A}_v$ graded projection operator which
quasi-permutes with $A$ so that ${\cal R}(F) \subset {\cal D}(A)$,
where ${\cal D}(A) = Domain (A)$, $~{\cal R}(A)=Range (A)$. Suppose
that $G$, $H$ and $J$ are the restrictions of $A$, $ ~ FB$ and $FD$
to ${\cal R}(F)$ respectively. Then $G$, $~H$ and $J$ are bounded
operators so that $~H$ and $J$ quasi-permute with $G$. Moreover,
$H^*$ and $J^*$ are the restrictions to ${\cal R}(F)$ of $B^*F^*$
and $C^*F^*$ respectively, where
\par $(1)$ $\mbox{}^j({\hat B}^*)\mbox{ }^k({\hat F}^*)=(-1)^{\kappa (j,k) +\eta (k)}
\mbox{ }^k{\hat F}\mbox{ }^j({\hat B}^*)$ and \par $(2)$
$\mbox{}^j({\hat D}^*) \mbox{ }^k({\hat F}^*) = (-1)^{\kappa (j,k)
+\eta (k)} \mbox{ }^k{\hat F} \mbox{ }^j
({\hat D}^*)$ \\
for each $j, k$, with $\kappa (j,k)=0$ for $j=k$ or $j=0$ or $k=0$,
$~\kappa (j,k) =1$ for $j\ne k\ge 1$, $~\eta (0)=0$, $~\eta (k)=1$
for each $k\ge 1$.}
\par {\bf Proof.}  Note 2.5 and Theorems 2.29, 2.44 and
Proposition 2.32 in \cite{ludunbnormopla12} and Definitions 1 imply
that in components the following formulas are satisfied:
$$(3)\quad \mbox{}_1^j{\hat E}(\delta _1)\mbox{ }_2^k{\hat E}(\delta _2) =
(-1)^{\kappa (j,k)} \mbox{ }_2^k{\hat E} (\delta _2)
\mbox{}_1^j{\hat E}(\delta _1)= \mbox{ }_2^j{\hat E} (\delta _2)
\mbox{}_1^k{\hat E}(\delta _1)$$ for each $j, k =0,1,2,...$, where
$\kappa (j,k)=0$ for $j=k$ or $j=0$ or $k=0$, $~\kappa (j,k) =1$ for
$j\ne k\ge 1$, \\ where $\theta ^j_k(x_j)$ is denoted by $x_j$ for
short, $\theta ^j_k: X_j\to X_k$ is an $\bf R$-linear topological
isomorphism of real normed spaces (see \S \S I.2.1 and I.2.73 in
\cite{ludopalglamb}). Suppose that $x, y \in {\cal R}(F)$, hence $x,
y \in {\cal D}(B) = {\cal D} (B^*)$, since ${\cal R}(F)\subset {\cal
D}(A)\subset {\cal D}(B)$. Therefore
\par $(4)$ $<FBx;y> = <Bx;y> = < x;B^*y> =<x;B^*y>=<x;FB^*y>$ and
\par $(5)$
$<\mbox{}^j{\hat F}\mbox{ }^k{\hat B}x_k;y_j> = <\mbox{}^k{\hat
B}x_k;y_j>i_j^* = (-1)^{\kappa (j,k)+\eta (k)} < x_k;\mbox{}^j{\hat
F}\mbox{ }^k({\hat B}^*)y_j>$.
\par If $L = FB^*|_{{\cal R}(F)}$, then $H^*=L$ and $H=L^*$ by
Formula $(4)$. The operator $L^*$ is closed, consequently, $H$ is
closed and ${\cal D}(H)\supset {\cal R}(F)$. In view of the closed
graph theorem for $\bf R$-linear operators the operator $H$ is
bounded 1.8.6 \cite{kadring}. This implies that the operator $G$ is
also bounded, since the operator $A$ is normal and hence closed so
that ${\cal R}(F)\subset {\cal D}(A)$. In view of Theorems 2.27,
2.29 an 2.44 in \cite{ludunbnormopla12} the operator $A$ has an
${\cal A}_v$ graded projection valued measure. Take now $x\in {\cal
R}(F)$, hence $Ax\in {\cal R}(F)\subset {\cal D}(A)\subset {\cal
D}(B)$, since
$$\mbox{}^j{\hat F} \mbox{ }^k{\hat F} = (-1)^{\kappa (j,k)} \mbox{
}^k{\hat F} \mbox{ }^j{\hat F}\mbox{ and }{\cal D}(F)={\cal
D}(F)_0i_0\oplus ... \oplus {\cal D}(F)_m i_m\oplus ...$$ for each
$j, k$ and
$$A=\int_{{\cal A}_v} F(dt).t$$ so that $\mbox{}^j{\hat F}\mbox{
}^k{\hat A}\subseteq (-1)^{\kappa (j,k)} \mbox{ }^k{\hat A}\mbox{
}^j{\hat F}$ for each $j, k$. Symmetric proof is for $A$ and $C$
instead of $A$ and $B$. The operators $B^*B$ and $C^*C$ belong to
the family $alg _{{\cal A}_v} (I,A,A^*,B,B^*,C,C^*)$.
\par In view of Theorem I.3.23 \cite{ludopalglamb} the spectra of $B^*B=\int_{-\infty
}^{\infty } \mbox{ }_{B^*B}F(dt).t^2 $ and $D^*D=\int_{-\infty
}^{\infty } \mbox{ }_{D^*D}F(dt).t^2 $ are real so that
$\mbox{}_{B^*B}F$ and $\mbox{}_{D^*D}F$ are ${\cal A}_v$ graded
projection valued measures for $B^*B$ and $D^*D$ respectively on
${\cal B}({\bf R})\subset {\cal B}({\cal A}_r)$.  Then from Formulas
$(2,4)$ and 1$(1, P1-P4)$ we deduce that
\par $(7)$ $(\mbox{}^j{\hat F}\mbox{ }^k{\hat B})\mbox{ }^s{\hat A}
x_s = (\mbox{}^j{\hat F} \mbox{ }^k{\hat B}) \sum_{p,q: ~i_pi_q=i_s}
[ \mbox{}^p{\hat D} \mbox{ }^q{\hat B}+ (-1)^{\kappa (p,q)} \mbox{
}^q{\hat D} \mbox{}^p
{\hat B}]$ \\
$= \sum_{p,q: ~i_pi_q=i_s} [ (\mbox{}^j{\hat F}\mbox{ }^k{\hat
B})(\mbox{}^p{\hat D} \mbox{ }^q{\hat B})+(-1)^{\kappa (p,q)}
(\mbox{}^j{\hat F}\mbox{ }^k{\hat B})(\mbox{}^q{\hat D} \mbox{
}^p{\hat B})]$\par $=(-1)^{\kappa (s,l)}\mbox{ }^s{\hat
A}(\mbox{}^j{\hat F} \mbox{ }^k{\hat B}x_s)$, \\
since the set theoretic composition of operators is associative:
$(FB)(DB)=F((BD)B)$, where $l$ is such that $i_ji_k\in {\bf R}i_l$.
Thus $H$ and analogously $J$ quasi-permute with $G$, since the
family $alg_{{\cal A}_v} (I,A,A^*, \mbox{ }_1A,\mbox{ }_1A^*, \mbox{
}_2A,\mbox{ }_2A)$ is quasi-commutative. From Formulas $(5,6)$ we
infer Equalities $(1,2)$.
\par {\bf 3. Notation.} Suppose that $a, b\in {\cal A}_r$. If
$b_j\ge a_j$ for each $j=0,1,2,...,2^r-1$, this fact will be denoted
by $b\succeq a$. Then ${\cal I}_{a,b} := \{ z\in {\cal A}_r:
b\succeq z \succeq a \} $.
\par {\bf 4. Lemma.} {\it Let operators $A$, $B$ and $D$ have
property $P$ and let $F$ be an ${\cal A}_v$ graded projection valued
measure for $A$, let also $b\succeq a\in {\cal A}_v$. Then ${\cal
R}(F({\cal I}_{a,b}))=:Y$ reduces both $B$ and $D$ and these
operators restricted to $Y$ are bounded and normal and they
quasi-permute with the restriction $A|_Y$.}
\par {\bf Proof.} Consider the pair of operators $A$ and $B$. Put
$\mbox{}_nF := F|_{{\cal I}_{-b(n),b(n)}}$ and $\mbox{}_nV={\cal
R}(F({\cal I}_{-b(n),b(n)}))$ with $b(n)_j = ni_j$ for every $n\in
\bf N$ and each $j=0,1,2,...,2^v-1$. Then $\mbox{}_nV\subset
\mbox{}_{n+1}V$ for each $n$. Therefore, an ${\cal A}_v$ vector
subspace $\bigcup_n \mbox{ }_nV=:V$ is dense in the Hilbert space
$X$ over the Cayley-Dickson algebra ${\cal A}_v$, consequently,
$\lim_n \mbox{ }_nF=I$ in the strong operator topology. Each
operator $\mbox{}_nA := A|_{\mbox{}_nV}$ is bounded and normal and
has the ${\cal A}_v$ graded projection valued measure on the Borel
$\sigma $-algebra ${\cal B}({\cal A}_v)$ of all Borel subsets in
${\cal A}_v$ so that $\mbox{}_nF= F|_{\mbox{}_nV}$ for each natural
number $n$. We consider the restriction $\mbox{}_nG :=
\mbox{}_nFB|_{\mbox{}_nV}$. It is known from Lemma 2, that each
operator $\mbox{}_nG$ is bounded and quasi-permutes with
$\mbox{}_nB$ so that $$(1)\quad \mbox{}_{n^G}\mbox{}^j{\hat F}
\mbox{ }^k_{n^B}{\hat F} = (-1)^{\kappa (j,k)} \mbox{ }^k_{n^B}{\hat
F}\mbox{ }^j_{n^G}{\hat F}$$ for each $j, k$, consequently,
$$(2)\quad \mbox{ }_n^s{\hat F}(\delta ) (\mbox{}^j_{n^B}{\hat F}(\delta _1)
\mbox{ }_n^k{\hat F}(\delta _2) x)=(-1)^{\kappa (j,k)} \mbox{
}_n^s{\hat F}(\delta ) (\mbox{ }_n^k{\hat F}(\delta _2) \mbox{
}^j_{n^B}{\hat F}(\delta _1)x)$$ for each $x\in \mbox{}_nV_0$ and
$\delta , \delta _1, \delta _2\in {\cal B}({\cal A}_v)$, where
$\mbox{}_{n^G}F$ and $\mbox{ }_{n^B}F$ denote ${\cal A}_v$ graded
projection valued measures for the operators $\mbox{}_nG$ and
$\mbox{}_nB$ correspondingly.
\par Let now $y\in {\cal D}(A)_0$ and $\delta \in {\cal B}({\cal
A}_v)$ be fixed, hence $$(3)\quad \lim_n \mbox{ }_n^s{\hat F}(\delta
) (\mbox{ }_n^k{\hat F}(\delta _2) \mbox{}^j_{n^B}{\hat F}(\delta
_1) x)=\lim_n ~ \mbox{ }^s{\hat F}(\delta ) (\mbox{ }_n^k{\hat
F}(\delta _2)\mbox{ }^j_{n^B}{\hat F}(\delta _1) x)$$  $$ =\pm
(-1)^{\psi (s,k,j)} \mbox{ }^l{\hat F}(\delta \cap \delta _2) \mbox{
}^j_{n^B}{\hat F}(\delta _1) x ,$$ where $i_si_k=\pm i_l$, $ ~ \psi
(s,j,k)\in \{ 0 , 1 \} $ is an integer so that
$i_s(i_ji_k)=(-1)^{\psi (s,j,k)} (i_si_j)i_k$. If a vector $x\in
\bigcup_n \mbox{ }_nV_0$ is given, then there exists a natural
number $m$ such that $$(4)\quad \mbox{ }_n^s{\hat F}(\delta )
(\mbox{ }^j_{n^B}{\hat F}(\delta _1) \mbox{ }_n^k{\hat F}(\delta _2)
x)= \mbox{}^s{\hat F}(\delta ) (\mbox{}^j_{n^B}{\hat F}(\delta _1)
\mbox{ }_n^k{\hat F}(\delta _2) x)$$ for each $n>m$, consequently,
$$(5)\quad \lim_n \mbox{ }_n^s{\hat F}(\delta ) (\mbox{}^j_{n^B}{\hat F}
(\delta _1) \mbox{ }_n^k{\hat F}(\delta _2) x) = \mbox{}^s{\hat
F}(\delta )(\mbox{}^j_B{\hat F}(\delta _1) \mbox{ }^k{\hat F}(\delta
_2) x),$$
 where $F({\cal A}_v)=I$, $~I$ denotes the unit operator.
 From Formulas $(2-5)$ and the
 inclusions $\bigcup_n \mbox{ }_nV =:V \subset {\cal D}(A)\subset {\cal
 D}(B)$ it follows, that
$$(6)\quad \mbox{}^j{\hat B}(\mbox{}^k{\hat F}(\delta ) x_si_s)=
(-1)^{\xi (j,k,s)} \mbox{ }^k{\hat F}(\delta )(\mbox{}^j{\hat
B}x_si_s)$$ for each $x_si_s\in V$ and $j,k, s=0,1,2,...$, where
$\xi (j,k,s)\in \{ 0, 1 \} $ is such integer number that
$i_j(i_ki_s)=(-1)^{\xi (j,k,s)}i_k(i_ji_s)$. From the formula
$i_j(i_ki_s)+i_k(i_ji_s)=2i_s Re(i_ji_k)$ we get $(-1)^{\xi
(j,k,s)}=(-1)^{\kappa (j,k)}$ for each $j$, $k$ and $s$, since an
algebra $alg_{\bf R}(i_j,i_k,i_s)$ over $\bf R$ generated by $i_j$,
$i_k$ and $i_s$ has an embedding into the octonion algebra which is
alternative \cite{baez} (see also Formulas 4.2.4$(7,8)$ in
\cite{ludancdnb}). Thus $BV\subset V$ and $\mbox{}_BFV\subset V$.
Then
$$(7)\quad \mbox{}^j_{n^H}{\hat F}(\delta _1)\mbox{ }_n^k{\hat F} (\delta )
= (-1)^{\kappa (j,k)} \mbox{ }_n^k{\hat F} (\delta )\mbox{
}^j_{n^H}{\hat F}(\delta _1) ,$$ that is $\mbox{}_nH^*$
quasi-permutes with $\mbox{}_nF$.
\par In view of Lemma 2 we have
$\mbox{}_nH^* = B^*|_{\mbox{}_nV}\mbox{}_nF({\cal I}_{-b(n),b(n)})$
and from the proof above we get $$(8)\quad \mbox{}_n^j{\hat F}
(\mbox{}^k_B({\hat F}^*) x_si_s) = (-1)^{\xi (j,k,s)} \mbox{
}^k_B({\hat F}^*) (\mbox{}_n^j{\hat F}x_si_s)$$ $$ = (-1)^{\xi
(j,k,s)}\mbox{ }^k_B({\hat F}^*) (\mbox{}^j{\hat I}x_si_s),$$
consequently, $B(\mbox{}_nV)\subset \mbox{}_nV$ and
$\mbox{}_BF(\mbox{}_nV)\subset \mbox{}_nV$. Consider decomposition
$x=y+z$ with $y\in \mbox{}_nV$ and $z\in \mbox{}_nV^{\perp }$, then
$x\in {\cal D}(B)$ is equivalent to $z\in {\cal D}(B)$. The latter
inclusion implies $z\in {\cal D}(B)\cap \mbox{ }_nV$, if
additionally $x\in \mbox{ }_nV$, then we get $<B^*y;z> = <y;Dz>=0$,
consequently, $Bz\in \mbox{ }_nV$ and this together with $(7)$ leads
to the inclusion $\mbox{}_nFB\subset B\mbox{ }_nF$, that is
$\mbox{}_n^j{\hat F}\mbox{ }^k{\hat B}\subset (-1)^{\kappa (j,k)}
\mbox{ }^k{\hat B} \mbox{ }_n^j{\hat F}$. For any $\bf R$-linear
spaces a sign in an inclusion does not play any role. Thus
$\mbox{}_nV$ reduces $B$ and $\mbox{}_nG$ into a normal operator
$\mbox{}_nQ=B|_{\mbox{}_nV}$.
\par Suppose that $\mbox{}_{n^G}F$ is the canonical ${\cal A}_v$ graded
projection valued measure for $\mbox{}_nG$ and $\mbox{}_BF$ is the
canonical ${\cal A}_v$ graded projection valued measure for $B$,
hence $\mbox{}_BF|_{\mbox{}_nV}=\mbox{}_{n^G}F$ for each $n\in \bf
N$. If $x\in \bigcup_n \mbox{ }_n V$, there exists a natural number
$m$ so that $$(9)\quad \mbox{}^j{\hat E}(\delta _1) (\mbox{}^k{\hat
F}(\delta _2) x_si_s)= \mbox{}^j_{n^G}{\hat F}(\delta _1)
(\mbox{}_n^k{\hat F}(\delta _2) x_si_s)$$ $$= (-1)^{\xi
(j,k,s)}\mbox{ }_n^k{\hat F}(\delta _2)(\mbox{}^j_{n^G}{\hat
F}(\delta _1)x_si_s)
$$  for each Borel subsets $\delta _1$ and $\delta _2$ in ${\cal
A}_r$, since the restriction of $A$ to $\mbox{}_nV$ and $\mbox{}_nG$
quasi-permute for all $n$ in accordance with Lemma 2. On the other
hand, the ${\cal A}_v$ vector space $V$ is dense in $X$,
consequently, $\mbox{}_BF$ and $F$ quasi-permute: $$(10)\quad
\mbox{}^j_B{\hat F}(\delta _1) \mbox{ }^k{\hat F}(\delta _2)x_0 =
(-1)^{\kappa (j,k)} \mbox{ }^k{\hat F}(\delta _2)\mbox{ }^j_B{\hat
F}(\delta _1)x_0 $$  for each $j, k=0,1,2,...$ and $x_0\in X_0$.
\par If now $F$ is an ${\cal A}_v$ graded projection valued measure
described in this lemma, then Formula $(10)$ implies $$(11)\quad
\mbox{}^k{\hat F} \mbox{}^j{\hat B}\subseteq (-1)^{\kappa (j,k)}
\mbox{ }^j{\hat B} \mbox{ }^k{\hat F}$$ for each $j, k=0,1,2,...,
2^v-1$, consequently, ${\cal R}(F)$ reduces $B$ and $B|_{{\cal
R}(F)}$ is a normal operator with ${\cal R}(F)\subset {\cal D}(B)$,
since ${\cal R}(F)\subset {\cal D}(A)\subset {\cal D}(B)$. This
restriction $B|_{{\cal R}(F)}$ is bounded by the closed graph
theorem 1.8.6 \cite{kadring}. Moreover, the restrictions of $A$ and
$B$ to ${\cal R}(F)$ quasi-permute. Analogous proof is valid for the
pair $A$ and $C$ instead of $A$ and $B$.

\par {\bf 5. Theorem.} {\it If operators $A$, $B$ and $D$
satisfy property $P$, then $B$ and $D$ quasi-permute so that
$$ (1)\quad \mbox{}^j{\hat B}\mbox{ }^k{\hat D} = (-1)^{\kappa (j,k)} \mbox{
}^k{\hat D}\mbox{ }^j{\hat B}$$ for each $j, k$. Moreover,
$$(2)\quad \mbox{}^l{\hat A} =\sum_{j,k; i_ji_k=i_l} (\mbox{}^j{\hat B}
\mbox{ }^k{\hat D} + (-1)^{\kappa (j,k)} \mbox{ }^k{\hat B}\mbox{
}^j{\hat D})$$ for each $l$.}
\par {\bf Proof.} Consider the canonical ${\cal A}_v$ graded
projection valued measure $E$ for a normal operator $A$ (see
Definition 1). Then we put $\mbox{}_nF := E({\cal I}_{a,b}))$ with
$a_j=-ni_j$ and $b_j=ni_j$ for each $j$. From Theorems 2.27, 2.29
and 2.44 in \cite{ludunbnormopla12} and \S 4 above we know that
$$(3)\quad Ax=\int_{{\cal A}_v} d\mbox{}_AE(t).t x\quad \forall x\in {\cal
D}(A)\mbox{  and}$$
$$(4)\quad Bx=\int_{{\cal A}_r} d\mbox{}_BE(t).t x\quad \forall x\in {\cal
D}(B)\mbox{  and}$$
$$(5)\quad Dx=\int_{{\cal A}_v} d\mbox{}_DE(t).t x\quad \forall x\in {\cal
D}(D),$$ where $\mbox{}_AE$, $\mbox{}_BE$ and $\mbox{}_DE$ denote
${\cal A}_v$ graded projection valued measures for $A$, $B$ and $D$
respectively. Then the condition $A=BD$ gives
$$(6)\quad Ax = \int_{{\cal A}_v} d\mbox{}_BE(t).t
\int_{{\cal A}_v} d\mbox{}_DE(u).u x .$$ To operators $A$, $B$ and
$D$ normal functions $h_A$, $h_B$ and $h_D$ correspond so that
$h_A=h_Bh_D$. On the other hand, to the operators $A^*A$ and $B^*B$
and $D^*D$ non-negative self-adjoint functions $|h_A|^2$, $|h_B|^2$
and $|h_D|^2$ correspond (see Proposition 2.32 in
\cite{ludunbnormopla12}). These operators $A$ and $B$ and $D$ are
normal so that they satisfy the identities
$A^*A=D^*B^*BD=D^*BB^*D=AA^*=BDD^*B^*=BD^*DB^*$ and $B^*B=B^*B$ and
$D^*D=DD^*$.
\par In view of Theorems 2.29, 2.44 and Proposition 2.32 and Remark
2.43 in \cite{ludunbnormopla12} to the ${\cal A}_v$ graded
projection operator $\mbox{}_AE(\delta )$ a homomorphism $\phi $ a
(real) characteristic function $\phi (\mbox{}_AE(\delta ))=\chi _u$
of a subset $u\subset \Lambda $ counterpose so that $\chi _u=\omega
(\chi _{\delta })$. Therefore, Theorem 2.23 and Lemma 2.21 in
\cite{ludunbnormopla12}, Formulas $(3-6)$ and Conditions $(P1-P4)$
imply that their projection operators satisfy the equality
\par $(7)$ $\mbox{}_BE(\delta _1)\mbox{}_DE(\delta _2) = \mbox{}_DE(\delta
_2) \mbox{}_BE(\delta _1)$ \\ for each Borel subsets $\delta _1$ and
$\delta _2$ in ${\cal A}_v$. In view of Lemma 4 ${\cal
R}(\mbox{}_nF)$ reduces $B$ and $D$ and the restrictions of these
operators to ${\cal R}(\mbox{}_nF)$ are bounded normal operators. On
the other hand, $\bigcup_{n=1}^{\infty } {\cal R}(\mbox{}_nF)$ is
dense in the Hilbert space $X$ over the Cayley-Dickson algebra
${\cal A}_v$. Therefore, we infer from Formulas $(3-7)$, that $
\mbox{}^jB$ and $\mbox{ }^kD$ satisfy Formulas $(1,2)$ for each $j,
k$, since
$$(8)\quad \mbox{}^j_B{\hat E}(\delta _1)\mbox{ }^k_D{\hat E}(\delta _2) =
(-1)^{\kappa (j,k)} \mbox{ }^k_D{\hat E}(\delta _2)\mbox{ }^j_B{\hat
E}(\delta _1)$$ for every Borel subsets $\delta _1$ and $\delta _2$
in ${\cal A}_v$ and for each $j, k$.
\par {\bf 6. Corollary.} {\it Suppose that operators $A$, $B$ and
$D$ are self-adjoint and satisfy property $(P)$. Then $BD=DB$.}
\par {\bf Proof.} This follows immediately from
Theorems 2.27, 2.29 and 2.44 in \cite{ludunbnormopla12} and Formulas
5$(1-3)$, since spectra of self-adjoint operators are contained in
the real field $\bf R$ and the latter is the center of the
Cayley-Dickson algebra ${\cal A}_v$ so that $t=t_0\in {\bf R}$ in
Formulas 5$(1,2)$, that is $j=k=0$ only.
\par {\bf 7. Lemma.} {\it Let operators $B$, $D$ and $A$ have property
$P$, let also $B=T_BU_B$, $~D=T_DU_D$ and $A=TU$ be their canonical
decompositions with positive self-adjoint operators $T_B$, $T_D$ and
$T$ and unitary operators $U_B$, $U_D$ and $U$ respectively. Then
$T_BT_D=T_DT_B=T$ and $U_BU_D=U$ so that $\mbox{}^jU_B\mbox{ }^kU_D=
(-1)^{\kappa (j,k)} \mbox{ }^jU_B\mbox{ }^kU_D$ for each $j, k$,
moreover, $T_BU_D=U_DT_B$ and $T_DU_B=U_BT_D$.}
\par {\bf Proof.} The decompositions in the conditions of this lemma
are particular cases of that of Theorem I.3.37 \cite{ludopalglamb}.
Consider the canonical ${\cal A}_v$ graded resolutions of the
identity $E^B$ and $E^D$ of operators $B$ and $D$ respectively. In
view of Theorem 5
$$\mbox{}^jE^B(\delta _1) \mbox{ }^kE^D(\delta _2) =
(-1)^{\kappa (j,k)} \mbox{ }^kE^D(\delta _2) \mbox{}^jE^B(\delta
_1)$$ for every Borel subsets $\delta _1$ and $\delta _2$ in ${\cal
A}_r$ and each $j, k$. We put $F(dw,dz) = E^B(dw)E^D(dz)$, hence
$F(dw,dz)$ is a $2^{v+1}$ parameter ${\cal A}_v$ graded resolution
of the identity so that $F_{i_k}(\delta _1,\delta _2)x_k =
E^B(\delta _1)(E^D)_{i_k}(\delta _2)$ for each vector $x_k\in X_k$
and every $k$ and we put
$$G:=\int_{{\cal A}_v^2} dF(w,z).wz ,$$  where $dF(w,z)$ is another
notation of $F(dw,dz)$, $~w, z\in {\cal A}_v$ (see also \S I.2.58
\cite{ludopalglamb}). This operator $G$ is normal, since the
quaternion skew field is associative and the octonion algebra is
alternative and $(wz)(wz)^*=|wz|^2=|w|^2|z|^2$ for each $w, z\in
{\cal A}_v$ with $2\le v\le 3$. Then we get $$B=\int_{{\cal A}_v^2}
dF(w,z).w=\int_{{\cal A}_v} dE^B(w).w\mbox{  and}$$
$$D=\int_{{\cal A}_v^2} dF(w,z).z=\int_{{\cal A}_v} dE^D(z).z
\mbox{, consequently,}$$ $A=BD$ and $\mbox{}^jB\mbox{ }^kD =
(-1)^{\kappa (j,k)}\mbox{ }^kD\mbox{}^jB$ for each $j, k$, and hence
$$\sum_{j, k: ~i_ji_k=i_l} [\mbox{}^jB\mbox{ }^kD + (-1)^{\kappa (j,k)}
\mbox{}^kB\mbox{ }^jD]\subseteq \mbox{ }^lG$$ for every $l$.
Therefore, $A=G$, since a normal operator is maximal. \par Then one
can consider the function $u(w,z):= wz/|wz|$ for $wz\ne 0$, while
$u(w,z)=1$ if $wz=0$, where $w, z \in {\cal A}_v$. The operator
$$U:=\int_{{\cal A}_v^2} dF(w,z).u(w,z)$$ is unitary, since
$|u(w,z)|=1$ for each $w$ and $z$, the operator
$$T := \int_{{\cal A}_v^2}  dF(w,z).|wz|$$ is positive and self-adjoint,
since $$<xT;x> := \int_{{\cal A}_v^2} <xdF(w,z).|wz| ;x>\ge 0$$ for
each $x\in {\cal D}(T)$ (see Proposition 2.35
\cite{ludunbnormopla12}). On the other hand,
$u(w,z)|wz|=|wz|u(w,z)=wz$, since the algebra ${\cal A}_v$ is
alternative for $v\le 3$, hence $TU=UT=G=A$ by Theorem 2.44
\cite{ludunbnormopla12}. Moreover, we deduce from Theorem 2.44
\cite{ludunbnormopla12} that the operators
$$U_B := \int_{{\cal A}_v^2} dF(w,z).u(w)\mbox{  and}$$
$$U_D := \int_{{\cal A}_v^2} dF(w,z).u(z)$$ are unitary and the
the operators
$$T_B := \int_{{\cal A}_v^2} dF(w,z).|w| \mbox{  and}$$
$$T_D := \int_{{\cal A}_r^2} dF(w,z).|z| $$ are positive and self-adjoint,
where $u(w) := w/|w|$ if $w\ne 0$, also $u(w)=1$ if $w=0$. Since
$|w||z|=|wz|$ for each $w$ and $z\in {\cal A}_v$ with $v\le 3$, the
inclusion follows $$T_BT_D\subseteq \int_{{\cal A}_v^2}
dF(w,z).|wz|=T.$$ The functions $u(w)$ and $u(z)$ are bounded and
$u(w)u(z)=u(z)u(w)=u(w,z)$ on ${\cal A}_v^2$, consequently,
$$U_BU_D=\int_{{\cal A}_v^2} dF(w,z).u(w,z)  =U\mbox{  so that}$$
$\mbox{}^jU_B\mbox{ }^kU_D= (-1)^{\kappa (j,k)} \mbox{ }^kU_B\mbox{
}^jU_D$ for each $j, k$. This implies that $A=UT=U_BU_DT=
(U_BT_B)(U_DT_D)=(U_BT_B)(T_DU_D)$, consequently, $U_DTU_D^*=
T_BT_D$. This means that the operators $T$ and $T_BT_D$ are
unitarily equivalent, hence the operator product $T_BT_D$ is
self-adjoint. A self-adjoint operator is maximal, consequently, $T
=T_BT_D$ and similarly $T=T_DT_B$. The real field $\bf R$ is the
center of the Cayley-Dickson algebra ${\cal A}_v$ for each $v\ge 2$,
the real and complex fields are commutative, hence
$$T_BU_D=\int_{{\cal A}_v^2}  dF  (w,z).(|w| u(z))= U_DT_B\mbox{
and}$$
$$T_DU_B=\int_{{\cal A}_v^2}  dF  (w,z).(|z| u(w))= U_BT_D.$$
\par {\bf 8. Notation.} Let $\Omega $ denote the set of all
$n$-tuples $x=(x_1,...,x_m,x_{m+1},...,x_n)$ such that $x_1,...,x_m$
are non-negative integers, while $x_{m+1},...,x_n$ are non-negative
real numbers with $\sum_{j=1}^n x_j> 0$. Relative to the addition
$x+y=(x_1+y_1,...,x_n+y_n)$ this set $\Omega $ forms a semi-group.
\par {\bf 9. Theorem.} {\it Suppose that $ \{ B^x: ~ x \in \Omega \}
$ is a weakly continuous semi-group of normal operators, that is
satisfying the following conditions:
\par $(1)$ $B^x$ is a normal operator acting on a Hilbert space
$X$ over the Cayley-Dickson algebra ${\cal A}_v$ for each element
$x\in \Omega $;
\par $(2)$ $B^xB^y = B^{x+y}$ for each $x, y \in \Omega $;
\par $(3)$ the ${\cal A}_v$ valued scalar product $<B^xf;g>$
is continuous in $x\in \Omega $ for each marked $f, g \in {\cal D}
:= \bigcap_{x\in \Omega } {\cal D}(B^x)$;
\par $(4)$ a family $alg_{{\cal A}_v} \{ I, B^x, (B^x)^*: ~ x \in
\Omega \} $ is over the algebra ${\cal A}_v$ with $~2\le v\le 3$.
Then a unique $2n$-parameter ${\cal A}_v$ graded resolution $ \{
\mbox{}_{(a_1,...,a_n;b_1,...,b_n)}{\hat F}: ~ a, b \in \Omega \} $
of the identity exists so that $\mbox{}_{(a,b)}{\hat F}=0$ if a
negative coordinate $a_k<0$ exists for some $k=1,...,n$, moreover,
$$(5)\quad B^x = \int_{{\bf R}^{2n}}
d\mbox{ }_{(a,b)}{\hat F}.\{ a^x ~ \exp [x_1M_1(b_1)b_1]... \exp
[x_nM_n(b_n)b_n] \} ,$$ where
$$a^x = \prod_{k=1}^n a_k^{x_k},$$
$M_s: {\bf R}^n\to {\cal S}_v := \{ z\in {\cal A}_v: ~ |z|=1, Re
(z)=0 \} $ is a Borel function for each $s$, $~a=(a_1,...,a_n)$.}
\par {\bf Proof.} In view of Lemma 5 each operator $B^x$ has the
decomposition $B^x=T^xU^x=U^xT^x$ with a positive self-adjoint
operator $T^x$ and a unitary operator $U^x$. Since $ \{ B^x: ~ x \in
\Omega \} $ is a semi-group, the relations $T^xT^y=T^{x+y}$ and
$U^xU^y=U^{x+y}$ are valid for each elements $x, y \in \Omega $.
That is, $ \{ T^x: ~ x \in \Omega \} $ and $ \{ U^x: ~ x \in \Omega
\} $ are semi-groups of positive self-adjoint operators and unitary
operators correspondingly.
\par If $y^s=(0,...,y^s_{m+1},...,y^s_n)\in \Omega $ are elements
of the semi-group $\Omega $ such that $y^1=\frac{y^2+y^3}{2}$,
$~s=1,2,3$, $ ~ f$ is a vector in a domain $\cal D$, then $$\|
B^{y^1}f\| ^2 = <B^{y^1}f,B^{y^1}f>
=<B^{y^2/2}B^{y^3/2}f;B^{y^2/2}B^{y^3/2}f>$$
$$=<(B^{y^2/2})^*B^{y^2/2}f; (B^{y^3/2})^*B^{y^3/2}f> \le \|
(B^{y^2/2})^*B^{y^2/2}f \|  \| (B^{y^3/2})^*B^{y^3/2}f \| $$ by
Cauchy-Schwartz' inequality I.2.4$(1)$ \cite{ludopalglamb}. On the
other hand, $$\| (B^{y^2/2})^*B^{y^2/2}f \| ^2 = <
(B^{y^2/2})^*B^{y^2/2}f;(B^{y^2/2})^*B^{y^2/2}f>=
<B^{y^2}f;B^{y^2}f>=\| B^{y^2}f \| ^2 ,$$ since the semi-group $\{
B^x: ~ x\in \Omega \} $ is commutative and an operator $B^x$ is
normal for each $x\in \Omega $. Thus the inequality $$\| B^{y^1}f\|
\le \| B^{y^2}f \| \| B^{y^3}f \| $$ follows. This implies that the
function $q(y) := \| B^yf \| $ is convex and bounded in the variable
$y_p$ in any bounded segment $[\alpha ,\beta ]\subset (0,\infty )$,
when other variables $y_q$ with $q\ne p$ are zero, $~p=m+1,...,n$,
since the exponential $e^t$ and the natural logarithmic functions
$\ln (t) $ are convex and bounded on each segment $[\gamma , \delta
]\subset (0,\infty )$ and $\ln q(y^1)\le \ln q(y^2)+\ln q(y^3)$.
\par Evidently, a commutative group $\hat \Omega $ exists for the semi-group
$\Omega $ such that $\Omega \subset {\hat \Omega }\subset {\bf R}^n$
and the function $q(y)$ can be extended on $\hat \Omega $ so that
$q(0)= \| f \| $ and $q(-y)=q(y)$ for $y\in \Omega $. If $q$ is
continuous on $\Omega $, its extension on $\hat \Omega $ can be
chosen continuous, since $\hat \Omega $ is a completely regular
topological space, i.e. $T_1$ and $T_{3.5}$ (see \cite{eng}).
\par If $\Omega $ is a group the function $q(y)$
is positive definite, that is by the definition for each $\lambda
_1,...,\lambda _k\in {\bf R}\oplus {\bf R}{\bf i}=:{\bf C}_{\bf i}$
and $y^1,...,y^k\in \Omega $ the inequality
$$\sum_{j,l} \lambda _j{\bar {\lambda }}_l q(y^j-y^l) \ge 0 $$ is valid,
but this inequality follows from the formula $$\sum_{j,l} \lambda
_j{\bar {\lambda }}_l q(y^j-y^l)= \| \sum_j\lambda _jB^{y^j}f \|
^2$$ and since $\| x \| \ge 0$ for each $x\in X$.
\par Particularly, for elements $x^k:=
(0,...,x_k,0,...,0)$ in the semi-group $\Omega $ the mapping
$<T^{x^k}f;f>$ is continuous in $x^k$ for each marked vector $f\in
{\cal D}$. Indeed, for $k=1,...,m$ this is evident, since $x^k\in
{\bf N}$ takes values in the discrete space in this case. If
$k=m+1,...,n$ one can use the formula $<T^{x^k}f;f>
=<B^{x^k/2}f,B^{x^k/2}f>= \| B^{x^k/2}f \| ^2$ which implies that
$<T^{x^k}f;f>$ is a bounded convex function of $x^k$ in every finite
interval $[\alpha ,\beta ]\subset (0,\infty )$, when $f\in \cal D$
is a marked vector (see Theorem 2.29 and Formula 2.44$(5)$
\cite{ludunbnormopla12}).
\par Denote by $\mbox{}_{s,t_s}E$ the canonical ${\cal A}_v$ graded
resolution of the identity for $T^{e_s}$, where
$e_s=(0,...,0,1,0,...)$ denotes the basic vector with coordinate $1$
at $s$-th place and zeros otherwise, $t_s\in \bf R$. By the
conditions of this theorem operators $T^{e_s}$ and $T^{e_p}$ commute
for each $s$, $p=1,...,n$, since
\par $(6)$ $T^{e_s}T^{e_p}=T^{e_s+e_p}=T^{e_p}T^{e_s}$. \\ Due to Theorem
2.42 \cite{ludunbnormopla12} the equality
$$(7)\quad \mbox{}^j_{s,t_s}E\mbox{ }^k_{p,t_p}E=(-1)^{\kappa (j,k)}
\mbox{ }^k_{p,t_p}E\mbox{ }^j_{s,t_s}E$$ is satisfied for each $j,
k$ and every $s, p$, with $t_s, t_p\in \bf R$. This implies that
\par $(8)$ $\mbox{}_{(t_1,...,t_n)}E = \mbox{}_{1,t_1}E...\mbox{}_{n,t_n}E$
\\ is an $n$-parameter ${\cal A}_v$ graded resolution of the identity.
Each operator $T^{e_s}$ is positive, hence $\mbox{}_{s,t_s}E=0$ for
every $t_s<0$, consequently, $\mbox{}_{(t_1,...,t_n)}E=0$ if $t_s<0$
for some $s=1,...,n$.
\par We now consider the operators
$$(9)\quad A^px := \int_0^{\infty }...\int_0^{\infty }
d\mbox{ }_{(t_1,...,t_n)}E.(t_1^{p_1}...t_n^{p_n}x),$$ where
$p=(p_1,...,p_n)\in \Omega $, $~x\in X$ for which the integral
converges. We certainly have $$\int_0^{\infty }...\int_0^{\infty }
d\mbox{ }_{(t_1,...,t_n)}E.(t_1^{p_1}...t_n^{p_n}x)= \int_0^{\infty
}...\int_0^{\infty } (t_1^{p_1}...t_n^{p_n}) d\mbox{
}_{(t_1,...,t_n)}E.x,$$ since $t_j^{p_j}\in {\bf R}$ for each $j$
and $\mbox{ }_{(t_1,...,t_n)}E$ is a real linear operator. If
$p_s\in {\bf Z}/2$ for each $s$, then $T^p=T^{e_1p_1}...T^{e_np_n}
\subseteq A^p$, consequently, $T^p=A^p$, since a self-adjoint
operator is maximal.
\par Take a partition of the Euclidean space ${\bf R}^n$
into a countable family of bounded parallelepipeds $J_k =
\prod_{j=1}^n [a_j,b_j]$ so that they may intersect only by their
boundaries: $J_k\cap J_l = \partial J_k\cap \partial J_l$ for each
$k\ne l\in \bf N$, $~\bigcup_{k=1}^{\infty }J_k={\bf R}^n$. We put
$Y^k := {\cal R}({\hat {\bf E}}(J_k))$, where ${\hat {\bf E}}(\delta
)$ is the ${\cal A}_v$ graded spectral measure corresponding to
$\mbox{}_tE$, $~\delta \in {\cal B}({\bf R}^n)$, $~t\in {\bf R}^n$.
Then the restriction $B^x|_{Y^k}$ of $B^x$ to $Y^k$ is a bounded
self-adjoint operator. If $x, y\in \Omega $ are elements of the
semi-group so that $y_s\ge x_s$ and $y_s\in {\bf Z}/2$ for each
$s=1,...,n$, then ${\cal D}(T^y)\subseteq {\cal D}(T^x)$, since
$T^y=T^xT^{y-x}$. Therefore, $f\in {\cal D}(A^y)={\cal
D}(T^y)\subseteq {\cal D}(T^x)$ for each $f\in Y^k$, consequently,
$Y^k\subset {\cal D}$ for each natural number $k\in \bf N$. \par If
$f\in Y^k\oplus Y^l$ and $g\in Y^l$, then $$\lim_{y\to x}
<(T^y-A^y)(f+g);(f+g)>=<(T^x-A^x)(f+g);(f+g)>=0,$$  since $T^y=A^y$
for each $y\in ({\bf Z}/2)^n\cap \Omega $ and the ${\cal A}_v$
valued scalar products $<T^xf;f>$ and $<A^xf;f>$ are continuous in
each component $x_s$ of $x$. In the same manner we get
$<(T^x-A^x)f;f>=0$ and $<(T^x-A^x)g;g>=0$, consequently,
$<(T^x-A^x)f;g>=0$. The ${\cal A}_v$ vector space
$\bigcup_{k=1}^{\infty }Y^k$ is dense in the Hilbert space $X$ over
the Cayley-Dickson algebra ${\cal A}_v$, hence
$T^xf^k=(A^x|_{Y^k})f^k=A^xf^k$ for each vector $f^k\in Y^k$. This
means that each $Y^k$ reduces the operator $T^x$ to $(A^x|_{Y^k})$,
consequently, $T^x=A^x$. From this it follows that the ${\cal A}_v$
valued scalar product $<T^xf;g>$ is continuous in $x\in \Omega $ for
each marked vectors $f\in {\cal D}$ and $g\in X$.
\par Consider the sub-semi-group $ \Omega _s := \{ x: ~
x=x^s := (0,...,0,x_s,0,...)\in \Omega \} $, where $s=1,...,n$, also
we suppose that $\mbox{}_s{\hat {\bf E}}( \{ 0 \} )=0$, where
$\mbox{}_s{\hat {\bf E}}( \delta
 )$ is the ${\cal A}_v$ graded projection valued measure corresponding to
 $\mbox{}_{s,t_s}E$, $~\delta \in {\cal B}({\bf R})$. This implies
 that the operator $T^{x^s}$ has not the zero eigenvalue. Take
 arbitrary marked vectors $f\in {\cal D}$ and $g\in {\cal
 D}(T^{y^s})$. Then using the triangle inequality we deduce that
 $$| <(U^{x^s}-U^{y^s})f;T^{y^s}g> | = |<(U^{x^s}-U^{y^s})T^{y^s}f;g>
 | = | <(U^{x^s}T^{x^s}-U^{y^s}T^{y^s})f;g> $$ $$+
 <U^{x^s}(T^{y^s}-T^{x^s})f;g> | \le | <(B^{x^s} - B^{y^s})f;g> |
 + \| (T^{y^s}-T^{x^s})f\| \| g \| . $$ But the limits are zero
$\lim_{x^s \to y^s } <(B^{x^s} - B^{y^s})f;g> =0$ due to
suppositions of this theorem and $\lim_{x^s \to y^s } \| (T^{y^s} -
T^{x^s})f\| =0$, since $T^x=A^x$ and $A^x$ has the integral
representation given by Formula $(9)$. Thus the limit $$\lim_{x^s
\to y^s } <(U^{x^s} - U^{y^s})f;h> =0$$ is zero for each $f\in \cal
D$ and $h\in {\cal R}(T^{y^s})$. On the other hand, $\cal D$ is
dense in $X$, since $\bigoplus_{k=1}^{\infty } Y^s$ is dense in $X$.
The family $U^x$ of unitary operators is norm bounded by the unit
$1$, consequently, $\lim_{x^s \to y^s } <(U^{x^s} - U^{y^s})f;h> =0$
for each $f, h\in X$ and hence the semi-group $ \{ U^{x^s}: ~ x^s\in
\Omega \} $ is weakly continuous. The semi-group $ \{ U^{x^s}: ~
x^s\in \Omega \} $ of unitary operators can be extended to a weakly
continuous group of unitary operators putting $U^{-x^s}=(U^{x^s})^*$
and $U^0=I$. This one-parameter commutative group of unitary
operators is also strongly continuous, since $$\| (U^{x^s}-U^{y^s})f
\| ^2 = <(U^{x^s}-U^{y^s})f;(U^{x^s}-U^{y^s})f> =$$
$$<(U^{x^s}-U^{y^s})^*(U^{x^s}-U^{y^s})f;f>=
<(2I-U^{x^s-y^s}-U^{y^s-x^s})f;f>$$  $$=
<(U^0-U^{x^s-y^s})f;f>+<(U^0-U^{y^s-x^s})f;f>.$$  \par In view of
Theorem I.3.28 \cite{ludopalglamb} there exists a unique ${\cal
A}_v$ graded projection valued measure $\mbox{}_s{\hat {\bf F}}$ so
that
$$(10)\quad <U(x^s)f;h> = \int^{\infty }_{-\infty }
<\mbox{}_s{\hat {\bf F}}(db_s).\exp (x_sM_s(b_s)b_s)f;h>$$ for each
$f, h \in {\cal D}(Q^s)$, where
$$(11)\quad <Q^sf,h> = \int^{\infty }_{-\infty } b_s <\mbox{}_s{\hat {\bf F}}(db_s)f;h>$$
for each $f, h \in {\cal D}(Q^s)$,
$$(12)\quad {\cal D}(Q^s)=\{ f: ~f\in X; ~ \| Q^sf \| ^2 =
\int^{\infty }_{-\infty }  <\mbox{}_s{\hat {\bf F}}(db_s).b_s^2f;f>
<\infty \} ,$$ $M_s(b_s)$ is a Borel function from $\bf R$ into the
purely imaginary unit sphere ${\cal S}_v := \{ z\in {\cal A}_v: ~
|z|=1, Re (z)=0 \} $. Then we put $\mbox{}_s{\hat {\bf E}}
(da_s,db_s)= \mbox{}_s{\hat {\bf E}}(da_s) \mbox{ }_s{\hat {\bf
F}}(db_s)$, where
$$\mbox{}_s^j{\hat {\bf E}}(\delta _1)
\mbox{ }_s^k{\hat {\bf F}}(\delta _2)=(-1)^{\kappa (j,k)} \mbox{
}_s^k{\hat {\bf E}}(\delta _1)\mbox{ }_s^j{\hat {\bf F}}(\delta
_2)$$ for each $j, k$ and Borel subsets $\delta _1, \delta _2\in
{\cal B}({\bf R})$. Then an operator $P^{x^s}$ exists prescribed by
the formula:
$$(13)\quad P^{x^s} = \int_{-\infty }^{\infty } \int_0^{\infty }
\mbox{}_s{\hat {\bf E}}(da_s,db_s).[a_s^{x_s} ~ \exp
(x_sM_s(b_s)b_s)].$$ This implies the inclusion $B^{x^s}\subseteq
P^{x^s}$, but a normal operator is maximal, consequently,
$B^{x^s}=P^{x^s}$ for each $s$ and $x^s\in \Omega $.
\par Suppose now that $\mbox{}_s{\hat {\bf E}}( \{ 0 \} )\ne 0$, consider the
null space $N^s := ker (B^{x^s})$ of $B^{x^s}$. To each ${\cal A}_v$
graded projection valued measure $\mbox{}_s{\hat {\bf E}}(\delta )$
associated with the family $alg_{{\cal A}_v} (I,B^x, (B^x)^*)$ a
real valued characteristic function in ${\cal N}(\Lambda ,{\bf R})$
corresponds, where $\delta \in {\cal B}({\bf R}^2)$, consequently,
$N^s$ is an ${\cal A}_v$ vector subspace in $X$. Let $X=N^s\oplus
K^s$, hence $K^s$ is an ${\cal A}_v$ vector space, since $N^s$ is
the ${\cal A}_v$ vector subspace of the ${\cal A}_v$ Hilbert space
$X$. Take the restrictions $B^{x^s}|_{N^s} =: B^{x^s,N}$ and
$B^{x^s}|_{K^s} =: B^{x^s,K}$ of $B^{x^s}$ to $N^s$ and $K^s$
correspondingly. This implies that the semi-group of normal
operators $\{ B^{x^s,K}: ~ x^s\in \Omega \} $ possesses  the
property that none of the operators $B^{x^s,K}$ has zero eigenvalue.
From Formula $(13)$ it follows, that there exists a two-parameter
resolution $\mbox{}_{s,K}{\hat {\bf E}}$ of the identity so that
$$(14)\quad B^{x^s,K} = \int_{-\infty }^{\infty } \int_0^{\infty }
 \mbox{}_{s,K}{\hat {\bf E}}(da_s,db_s).[a_s^{x_s} ~ \exp
(x_sM_s(b_s)b_s)].$$ Define an ${\cal A}_v$ graded projection value
measure $\mbox{}_{s,N}{\hat {\bf E}}$ so that
$$\int_{-\infty }^{a_s}\int_{-\infty }^{b_s} \mbox{}_{s,N}{\hat
{\bf E}}(dt_s,dq_s) = \mbox{}_{s,N;a_s,b_s}{\hat {\bf E}}$$ so that
$\mbox{}_{s,N;a_s,b_s}{\hat {\bf E}}=0$ for $a_s<0$ and
$\mbox{}_{s,N;a_s,b_s}{\hat {\bf E}}=I$ when $a_s\ge 0$.
 Since $B^{x^s,N}(N^s)= \{ 0 \} $, the integral representation
 follows:
$$(15)\quad B^{x^s,N} = \int_{-\infty }^{\infty } \int_0^{\infty }
\mbox{}_{s,N}{\hat {\bf E}}(da_s,db_s).[a_s^{x_s} ~ \exp
(x_sM_s(b_s)b_s)].$$ Now it is natural to put $\mbox{}_s{\hat {\bf
E}}(da_s,db_s) = \mbox{}_{s,N}{\hat {\bf E}}(da_s,db_s)\oplus
\mbox{}_{s,K}{\hat {\bf E}}(da_s,db_s)$ for an ${\cal A}_v$ graded
projection valued measure on $X$, that induces the formula:
$$(16)\quad B^{x^s} = \int_{-\infty }^{\infty } \int_0^{\infty }
\mbox{}_s{\hat {\bf E}}(da_s,db_s).[a_s^{x_s} ~ \exp
(x_sM_s(b_s)b_s)].$$ In accordance with Theorem 5 $$\mbox{}_s^j{\hat
{\bf E}}(\delta _1)\mbox{ }_q^k{\hat {\bf E}}(\delta _2) =
(-1)^{\kappa (j,k)} \mbox{ }_q^k{\hat {\bf E}}(\delta
_2)\mbox{}_s^j{\hat {\bf E}}(\delta _1)$$ for each $s, q=1,...,n$
and $j, k=0,1,...,2^v-1$ and every $\delta _1, \delta _2\in {\cal
B}({\bf R}^2)$, particularly, for $j=k=0$ i.e. $\mbox{}_s{\hat {\bf
E}}(\delta _1)$ and $\mbox{ }_q{\hat {\bf E}}(\delta _2)$ commute.
Then we put
$$\mbox{ }_{(a,b)}{\hat F}= \int_{-\infty }^{a_1}\int_{-\infty
}^{b_1}...\int_{-\infty }^{a_n}\int_{-\infty }^{b_n}\mbox{}_1{\hat
{\bf E}}(dt_1,dq_1)...\mbox{}_n{\hat {\bf E}}(dt_n,dq_n),$$ hence
$\mbox{ }_{(a,b)}{\hat F}$ is an ${\cal A}_v$ graded resolution of
the identity, for which $$d\mbox{ }_{(a,b)}{\hat F}.\{ a^x ~ \exp
[x_1M_1(b_1)b_1]... \exp [x_nM_n(b_n)b_n] \} =$$
$$\mbox{}_1{\hat {\bf E}}(da_1,db_1).\exp (x_1M_1(b_1)b_1)...\mbox{}_n{\hat
{\bf E}}(da_n,db_n).\exp (x_nM_n(b_n)b_n)a^x,$$ since the
semi-groups $ \{ B^x: ~ x \in \Omega \} $ and  $ \{ T^x: ~ x \in
\Omega \} $ and $ \{ U^x: ~ x \in \Omega \} $ are commutative, the
real field $\bf R$ is the center of the Cayley-Dickson algebra
${\cal A}_v$ for each $v\ge 2$, the fields ${\cal A}_0=\bf R$ and
${\cal A}_1=\bf C$ are commutative, $a_s\in \bf R$ and $x_s\in \bf
R$ for each $s=1,...,n$. For the operators
$$(17)\quad P^x = \int_{{\bf R}^{2n}} d\mbox{ }_{(a,b)}{\hat F}.\{ a^x ~ \exp [x_1M_1(b_1)b_1]...
\exp [x_nM_n(b_n)b_n] \} ,$$ where
$$a^x = \prod_{k=1}^n a_k^{x_k},$$
$M_s: {\bf R}^n\to {\cal S}_v := \{ z\in {\cal A}_v: ~ |z|=1, Re
(z)=0 \} $ is a Borel function for each $s$, the inclusion follows
$B^x\subseteq P^x$ for each $x\in \Omega $, since
$B^x=B^{x^1}...B^{x^n}$, where the operators $B^{x^1}$,...,$B^{x^n}$
pairwise commute. But a normal operator is maximal, consequently,
$B^x=P^x$ for each $x\in \Omega $. A uniqueness of the resolution
$\mbox{ }_{(a,b)}{\hat F}$ of the identity follows from uniqueness
of $\mbox{}_s{\hat {\bf F}}$ and $\mbox{}_s{\hat {\bf E}}$ for each
$s$.

\end{document}